\documentclass[twoside,11pt
]{amsart}%

  \usepackage{amsmath}
  \usepackage{amssymb}
  \usepackage{amsthm}
  \usepackage{mathrsfs}
  \usepackage{bbm}
  \usepackage{pifont}
  \usepackage{color}
  \usepackage{graphicx}

	\def\CC{{\ifmmode{\mathbbm{C}}\else{$\mathbbm{C}$}\fi}}
	\def\EE{{\ifmmode{\mathbbm{E}}\else{$\mathbbm{E}$}\fi}}
	\def\FF{{\ifmmode{\mathbbm{F}}\else{$\mathbbm{F}$}\fi}}
	\def\HH{{\ifmmode{\mathbbm{H}}\else{$\mathbbm{H}$}\fi}}
	\def\KK{{\ifmmode{\mathbbm{K}}\else{$\mathbbm{K}$}\fi}}
	\def\NN{{\ifmmode{\mathbbm{N}}\else{$\mathbbm{N}$}\fi}}
	\def\PP{{\ifmmode{\mathbbm{P}}\else{$\mathbbm{P}$}\fi}}
	\def\QQ{{\ifmmode{\mathbbm{Q}}\else{$\mathbbm{Q}$}\fi}}
	\def\RR{{\ifmmode{\mathbbm{R}}\else{$\mathbbm{R}$}\fi}}
	\def\TT{{\ifmmode{\mathbbm{T}}\else{$\mathbbm{T}$}\fi}}
	\def\UU{{\ifmmode{\mathbbm{U}}\else{$\mathbbm{U}$}\fi}}
	\def\ZZ{{\ifmmode{\mathbbm{Z}}\else{$\mathbbm{Z}$}\fi}}
	
	\def\A{{\ifmmode{\mathscr{A}}\else{$\mathscr{A}$}\fi}}
	\def\B{{\ifmmode{\mathscr{B}}\else{$\mathscr{B}$}\fi}}
	\def\C{{\ifmmode{\mathscr{C}}\else{$\mathscr{C}$}\fi}}
	\def\D{{\ifmmode{\mathscr{D}}\else{$\mathscr{D}$}\fi}}
	\def\E{{\ifmmode{\mathscr{E}}\else{$\mathscr{E}$}\fi}}
	\def\F{{\ifmmode{\mathscr{F}}\else{$\mathscr{F}$}\fi}}
	\def\G{{\ifmmode{\mathscr{G}}\else{$\mathscr{G}$}\fi}}
	\def\H{{\ifmmode{\mathscr{H}}\else{$\mathscr{H}$}\fi}}
	\def\I{{\ifmmode{\mathscr{I}}\else{$\mathscr{I}$}\fi}}
	\def\J{{\ifmmode{\mathscr{J}}\else{$\mathscr{J}$}\fi}}
	\def\K{{\ifmmode{\mathscr{K}}\else{$\mathscr{K}$}\fi}}
	\def\L{{\ifmmode{\mathscr{L}}\else{$\mathscr{L}$}\fi}}
	\def\M{{\ifmmode{\mathscr{M}}\else{$\mathscr{M}$}\fi}}
	\def\N{{\ifmmode{\mathscr{N}}\else{$\mathscr{N}$}\fi}}
	\def\O{{\ifmmode{\mathscr{O}}\else{$\mathscr{O}$}\fi}}
	\def\P{{\ifmmode{\mathscr{P}}\else{$\mathscr{P}$}\fi}}
	\def\Q{{\ifmmode{\mathscr{Q}}\else{$\mathscr{Q}$}\fi}}
	\def\R{{\ifmmode{\mathscr{R}}\else{$\mathscr{R}$}\fi}}
	\def\S{{\ifmmode{\mathscr{S}}\else{$\mathscr{S}$}\fi}}
	\def\T{{\ifmmode{\mathscr{T}}\else{$\mathscr{T}$}\fi}}
	\def\U{{\ifmmode{\mathscr{U}}\else{$\mathscr{U}$}\fi}}
	\def\V{{\ifmmode{\mathscr{V}}\else{$\mathscr{V}$}\fi}}
	\def\W{{\ifmmode{\mathscr{W}}\else{$\mathscr{W}$}\fi}}
	\def\X{{\ifmmode{\mathscr{X}}\else{$\mathscr{X}$}\fi}}
	\def\Y{{\ifmmode{\mathscr{Y}}\else{$\mathscr{Y}$}\fi}}
	\def\Z{{\ifmmode{\mathscr{Z}}\else{$\mathscr{Z}$}\fi}}

	\newtheoremstyle{slanted}
	{}
	{}
	{\slshape}
	{}
	{\bfseries}
	{.}
	{ }
	{}
	
	\theoremstyle{slanted}
	\newtheorem{theo}{Theorem}[section]
	\newtheorem{prop}[theo]{Proposition}
	
	\newtheorem{question}[theo]{Question}
	
	\newtheorem{lemma}[theo]{Lemma}
	\newtheorem{defi}[theo]{Definition}
	\newtheorem{corollary}[theo]{Corollary}

	\newtheorem{property}{Property}
	\def\ind#1{\mathbbmss{1}_{#1}}
	\def\egdef{:=}
	\def\Id{\mathop{\mbox{Id}}}
	\def\Aut{\mathop{\text{Aut}}}
	\def\build#1_#2^#3{\mathrel{\mathop{\kern 0pt#1}\limits_{#2}^{#3}}}
	\def\tend#1#2{\build\hbox to 12mm{\rightarrowfill}_{#1\rightarrow #2}^{}}

\title{An introduction to joinings in ergodic theory}
\author{Thierry de la Rue}
\address{Laboratoire de Math\'ematiques Rapha\"el Salem\\
	UMR 6085 CNRS -- Universit\'e de Rouen\\
	Avenue de l'Universit\'e\\
	B.P. 12\\
	F76801 Saint-\'Etienne-du-Rouvray Cedex}
\email{thierry.de-la-rue@univ-rouen.fr}

\begin{document}
\bibliographystyle{amsplain}

\maketitle

\begin{abstract}
Since their introduction by Furstenberg \cite{furst1}, joinings have proved
a very powerful tool in ergodic theory. We present here some aspects of the use
of joinings in the study of measurable dynamical systems,
emphasizing on
\begin{itemize}
\item the links between the existence of a non trivial common factor and the existence of
a joining which is not the product measure,
\item how joinings can be employed to provide elegant proofs of classical results,
\item how joinings are involved in important questions of ergodic theory, such as
pointwise convergence or Rohlin's multiple mixing problem.
\end{itemize}

\end{abstract}

\tableofcontents

\section{What are we talking about?}

\subsection{Dynamical systems and stationary processes}

We call here a \emph{dynamical system} any quadruple of the form $(X,\A,\mu,T)$, where 
$(X,\A,\mu)$ is a Lebesgue space (or equivalently a standard Borel space equipped with a 
probability measure $\mu$) and $T$ an automorphism of $(X,\A,\mu)$: $T$ is a one-to-one measurable
transformation of $X$ satisfying, for any measurable subset $A$ of $X$,
$$ \mu(T^{-1}A)\ =\ \mu(TA)\ =\ \mu(A). $$
Throughout this text, we will often use simply the symbol $T$ to designate the 
dynamical system $(X,\A,\mu,T)$. We will also often need a second dynamical system $S$, which will stand
for the quadruple $(Y,\B,\nu,S)$

Such dynamical systems, which are the objects of interest in measurable ergodic theory, can 
be considered from a rather probabilistic viewpoint by studying \emph{stationary processes}. Indeed,
any measurable map $\xi_0:\ X\to E$ gives rise to a stationary process $\xi=(\xi_k)_{k\in\ZZ}$ defined
by $\xi_k\egdef \xi_0\circ T^k$. (The set $E$ in which $\xi$ takes its values could be any standard Borel space;
most of the time, it will be a finite or countable alphabet.)
Conversely, to any $E$-valued stationary process $\xi=(\xi_k)_{k\in\ZZ}$, we can associate in a canonical
way a dynamical system: take $X_{\xi}\egdef E^{\ZZ}$, the sample space of the whole process, equipped
with its Borel $\sigma$-algebra $\A_{\xi}$ and the law $\mu_{\xi}$ of the process. Since
$\xi$ is stationary, $\mu_\xi$ is invariant by the shift $T_{\xi}$:
$(x_k)_{k\in\ZZ}\longmapsto(x'_k)_{k\in\ZZ}$ where $x'_k\egdef x_{k+1}$. 

\begin{defi}
In the dynamical system $(X,\A,\mu,T)$, the stationary process $\xi=(\xi_0\circ T^k)_{k\in\ZZ}$
is called a \emph{generating process} if the $\sigma$-algebra it generates is, modulo $\mu$, the whole $\sigma$-algebra $\A$.
\end{defi}

Observe that in the dynamical system $(X_{\xi},\A_{\xi},\mu_{\xi},T_{\xi})$, the process $\xi$ (obtained by 
looking at the coordinates) is always a generating process.

The two dynamical systems $T$ and $S$ are said to be \emph{isomorphic} if we can find
invariant subsets\footnote{In the sequel, we will not explicitely mention this restriction 
to subsets of full measure; it should be understood that every map is really defined on a subset of
probability~1.} $X_0\subset X$, $Y_0\subset Y$, with $\mu(X_0)=\nu(Y_0)=1$, and
a measurable one-to-one map $\psi:\ X_0\to Y_0$ with, for all $B\in\B$, $\mu(\psi^{-1}B)=\nu(B)$ and such that
$S\circ \psi=\psi\circ T$. If $\xi$ is a generating process for the dynamical system $T$, then $T$ is 
isomorphic to the dynamical system $T_\xi$ defined above. Recall the following well-known result of
ergodic theory, which shows that studying dynamical systems is essentially the same thing as studying
stationary processes taking values in a countable alphabet. (A nice proof can be found in \cite{kalik5}.)

\begin{theo}
\label{generating-process-theorem}
Every dynamical system admits a generating process taking its values in a countable
alphabet.
\end{theo}

\subsection{Products and factors: arithmetic of dynamical systems}

In his famous article initiating the theory of joinings \cite{furst1}, Furstenberg observes that 
a kind of arithmetic can be done with dynamical systems: Indeed, there are two natural
operations in ergodic theory which present some analogy with the integers.

First, if we are given two dynamical systems $(X,\A,\mu,T)$ and $(Y,\B,\nu,S)$, we can construct
their \emph{direct product} $(X\times Y, \A\otimes \B, \mu\otimes \nu, T\times S)$ where 
$T\times S:\ (x,y)\longmapsto(Tx,Sy)$. Note that the Cartesian product $X\times Y$ is here equipped with 
the direct product of the probability measures $\mu$ and $\nu$ , which is always $T\times S$-invariant. 
Acting on equivalence classes of isomorphic dynamical systems, the direct product operation  
is commutative, associative, and possesses a neutral element which is 
the trivial system reduced to one singleton. 
Translated in the context of stationary processes, if $\xi$ and $\zeta$ are respectively $E$-valued 
and $F$-valued stationary processes (not necessarily defined on the same probability space), 
their direct product $\xi\otimes\zeta$ is the $E\times F$-valued stationary process whose law is 
the direct product of the laws of $\xi$ and $\zeta$. That is to say, we are making the two
processes $\xi$ and $\zeta$ live together in an independant (hence stationary) way.

Second, we say that $(Y,\B,\nu,S)$ is a \emph{factor} of $(X,\A,\mu,T)$ if we can find
a measurable map $\pi:\ X\to Y$ satisfying
\begin{itemize}
\item $\pi(\mu)=\nu$ ;
\item $\pi\circ T=S\circ \pi$.
\end{itemize}
Such a map is called a \emph{homomorphism} of dynamical systems.
Heuristically, the existence of such a homomorphism means that we can see the system $S$ 
inside the system $T$ by looking at $\pi(x)$. 
To this factor is associated the \emph{factor $\sigma$-algebra} $\pi^{-1}(\B)$. This $\sigma$-algebra
has the property that if $A$ is measurable with respect to $\pi^{-1}(\B)$, then so are $T(A)$ and 
$T^{-1}(A)$. Conversely, each sub-$\sigma$-algebra $\F$ of $\A$ which is invariant by $T$ and $T^{-1}$
is a factor $\sigma$-algebra. Indeed, such a $\sigma$-algebra is always generated by some stationary
process $\xi$ living in the system $T$ (this is an adaptation of Theorem~\ref{generating-process-theorem}
to the case of a factor $\sigma$-algebra), and the dynamical system $T_\xi$ associated to $\xi$
is a factor  of $T$: Take $\pi(x):=(\ldots,\xi_{-1}(x),\xi_0(x),\xi_1(x),\ldots)$. 

Clearly, both systems $S$ and $T$ are factors of the direct product $T\times S$: These systems
are seen by looking at the two coordinates in the Cartesian product. The analogy with the
integers introduced by Furstenberg now gives two ways to define 
``$T$ and $S$ are relatively prime". 

\begin{property}
\label{no-common-factor-property}
$T$ and $S$ have no common factor except the trivial system reduced to one point.
\end{property}

\begin{property}
\label{disjointness-property}
Each time $T$ and $S$ appear as factors in some dynamical system, then their direct product
$T\times S$ also appears as a factor\footnote{%
To be exact, Property~\ref{disjointness-property}
should be stated in the following, stronger form: Each time $T$ and $S$ appear as factors in some 
dynamical system through the respective homomorphisms $\pi_T$ and $\pi_S$, $T\times S$ also appears 
as a factor through a homomorphism $\pi_{T\times S}$ such that $\pi_X\circ\pi_{T\times S}=\pi_T$ 
and $\pi_Y\circ\pi_{T\times S}=\pi_S$, where $\pi_X$ and $\pi_Y$ are the projections on the coordinates 
in the Cartesian product $X\times Y$.}.
\end{property}

The relationships between these two properties are exposed in Section~\ref{joinings-factors-section}.
Property~\ref{disjointness-property} has been called by Furstenberg the \emph{disjointness} property of
$T$ and $S$. The study of disjointness turns out to be quite enlightening, as it refers to \emph{all the
possible ways that two systems can be seen living together inside a third one}.  As explained in the 
next subsection, this is precisely the theory of joinings.

The large scale of applications of disjointness and joinings is already visible in Furstenberg's article
where the concepts of disjointness and joinings where introduced. Indeed, Furstenberg's initial motivations
ranged from the classification of dynamical systems (how disjointness can be used to characterized some
classes of dynamical systems) to a question in Diophantine approximation (which multiplicative semigroups of integers $S$ have the property that any real number is the limit of a sequence of rational numbers with denominators belonging to $S$?), passing through a filtering problem in probability theory (if $(\xi_n)$ and 
$(\zeta_n)$ are two real-valued stationary processes, can we recover $(\xi_n)$ from the knowledge of
$(\xi_n + \zeta_n)$?). Since Furstenberg's article, disjointness and joinings have been widely studied, and many other related notions have been introduced. In the present work, we shall mainly concentrate on some links between joinings and other ergodic properties of dynamical systems. For a more complete treatment of ergodic theory via joinings, we refer the readers to Eli Glasner's book \cite{glasn4}.

\subsection{Joinings}

\begin{defi}
Let $(X,\A,\mu,T)$ and $(Y,\B,\nu,S)$ be two dynamical systems. A \emph{joining} of $T$ and $S$
is a probability measure $\lambda$ on the Cartesian product $X\times Y$, whose marginals on $X$ and $Y$
are $\mu$ and $\nu$ respectively, and which is invariant by the product transformation $T\times S$.
\end{defi}

The set of joinings of $T$ and $S$, denoted by $J(T,S)$, is never empty since it always contains 
the product measure $\mu\otimes\nu$. Each $\lambda\in J(T,S)$ provides a new dynamical system 
$(X\times Y,\A\otimes\B,\lambda,T\times S)$, which will be denoted by $(T\times S)_\lambda$ 
to avoid confusion with the direct product $T\times S$. 
Both $T$ and $S$ are factors of $(T\times S)_\lambda$. Conversely, if $T$ and $S$ appear as factors in 
a third system $(Z,\C,\rho,R)$ through respectively the maps $\pi_X:\ Z\to X$ and $\pi_Y:\ Z\to Y$, we can 
consider the map $\pi:\ Z\to X\times Y$ defined by $\pi(z)\egdef(\pi_X(z),\pi_Y(z))$. Then the probability
measure $\lambda\egdef\pi(\rho)$ on $X\times Y$ satisfies the conditions to be a joining of $T$ and $S$, 
and the system $(T\times S)_\lambda$ is a factor of $R$. 

\subsection*{Topology on the set of joinings}

The set of joinings of $T$ and $S$ is equipped with the topology defined by
\begin{equation}
\label{def-topology}
\lambda_n\ \tend{n}{\infty}\ \lambda\ 
\Longleftrightarrow\ \forall A\in\A,\ \forall B\in\B,\ \lambda_n(A\times B)\ \tend{n}{\infty}\ \lambda(A\times B).
\end{equation}
This topology is metrizable: An example of distance defining it is given by
$$ d(\lambda,\lambda')\ \egdef\ \sum_{m,n\ge 0}\frac{1}{2^{m+n}}|\lambda(A_m\times B_n)-\lambda'(A_m\times B_n)|, $$
where $(A_m)_{m\ge0}$ and $(B_n)_{n\ge0}$ are countable algebras generating respectively the $\sigma$-algebras $\A$ 
and $\B$. With this topology, the set of joinings $J(T,S)$ is turned into a \emph{compact metrizable topological space}. It is interesting to observe that, when $X$ and $Y$ are compact metric spaces, this topology coincides with the weak topology on the set of joinings.

All these definitions can be extended in an obvious way to the notion of \emph{joining of a finite or countable 
family of dynamical systems} $(T_i)_{i\in I}$. If all the $T_i$'s are copies of a single dynamical system $T$, we
speak of \emph{self-joinings} of $T$. The set of self-joinings of order $k$ of $T$ (joinings of $k$ copies of $T$)
is denoted by $J_k(T)$, and we simply write $J(T)$ for $J_2(T)$. 

\subsection*{Ergodic joinings}

We denote by $J_e(T,S)$ the set of \emph{ergodic} joinings of $T$ and $S$. Since every factor of an ergodic system is still ergodic, a necessary condition for $J_e(T,S)$ not to be empty is that both $T$ and $S$ be ergodic. The following proposition shows that the converse is true.

\begin{prop}
  Given two ergodic dynamical systems $T$ and $S$, there exists at least one ergodic joining of $T$ and $S$. 
\end{prop}

\begin{proof}
We start from the only joining whose existence is known, $\lambda\egdef \mu\otimes\nu$. This joining may not be 
ergodic, but in this case we consider its \emph{decomposition in ergodic components}: There exists a 
probability measure $P$ on the set of all $T\times S$-invariant ergodic probability measures such that
$$ \lambda\ =\ \int \lambda_\omega\, dP(\omega). $$
Denoting by $\lambda^1_\omega$ (respectively $\lambda^2_\omega$) the marginals of the ergodic component
$\lambda_\omega$ on $X$ (respectively $Y$), we get
$$ \mu\ =\ \int \lambda^1_\omega\, dP(\omega),$$
and
$$ \nu\ =\ \int \lambda^2_\omega\, dP(\omega).$$
Since $\mu$ and $\nu$ are ergodic, this implies that for $P$-almost every $\omega$, 
$\lambda^1_\omega=\mu$ and $\lambda^2_\omega=\nu$. In other words, 
for $P$-almost every $\omega$, $\lambda_\omega$ is an ergodic joining of
$T$ and $S$. 
\end{proof}

Note that this proof implies a stronger result: When $T$ and $S$ are ergodic,
every joining of $T$ and $S$ is a convex combination of ergodic joinings. Since an
ergodic measure cannot be written as a non trivial convex combination of invariant
measures, we get the following.
\begin{prop}
\label{characterization-of-ergodic-joinings}
If $T$ and $S$ are ergodic, the set of ergodic joinings of $T$ and $S$ is the set 
of extremal points in the compact convex set $J(T,S)$.
\end{prop}

\section{From disjointness to isomorphism}

In this section, we present the two extreme cases concerning the joinings of
two dynamical systems $T$ and $S$. The first one is when the set of joinings $J(T,S)$ is reduced
to the singleton $\{\mu\otimes\nu\}$. Heuristically, this means that $T$ and $S$ have nothing 
in common which could lead to a non-trivial joining (what they have in common when they are not
disjoint is developped in Sections~\ref{dlcf-subsection} and~\ref{T-factor-subsection}). 
To the opposite, if $T$ and $S$ are isomorphic, we can construct
very special joinings of $T$ and $S$, namely the joinings supported on graphs of isomorphisms.

\subsection{Product measure and disjointness}

\begin{defi}
\label{disjointness-definition}
The dynamical systems $(X,\A,\mu,T)$ and $(Y,\B,\nu,S)$ are \emph{disjoint} if
the product measure $\mu\otimes\nu$ is the only joining of $T$ and $S$. We write in this
case $T\bot S$.
\end{defi}

We leave as an exercise to the reader the verification of the equivalence of this definition
with the strong form of Property~\ref{disjointness-property}. By the way,
why was it necessary to state this strong form? (Hint: There exists a non-trivial dynamical system  $T$
which is isomorphic to the direct product $T\times T$.)

We now give the most elementary example of disjointness: The identity is disjoint from 
any ergodic dynamical system.

\begin{prop}
If $T$ is the identity map on $X$ and $S$ is any ergodic dynamical system, then $T\bot S$.
\end{prop}

\begin{proof}
Let $\lambda$ be any joining of $T=\Id$ and $S$. For any $A\in\A$, $B\in\B$, invariance of $\lambda$
by $T\times S$ gives 
\begin{eqnarray*}
\forall n\ge 1,\quad \lambda(A\times B) & = & \lambda(A\times S^{-n}B)\\
  & = & \frac{1}{n}\sum_{k=0}^{n-1} \lambda(A\times S^{-k}B)\\
  & = & \int_{A\times Y} \frac{1}{n}\sum_{k=0}^{n-1} \ind{B}(S^k y)\,d\lambda(x,y).
\end{eqnarray*}
Since $S$ is ergodic, $\frac{1}{n}\sum_{k=0}^{n-1} \ind{B}(S^k y)$ converges to $\nu(B)$ $\nu$-almost everywhere,
hence $\lambda$-almost everywhere. This implies 
$$ \lambda(A\times B)=\mu(A)\nu(B), $$ 
hence $\lambda=\mu\otimes\nu$.
\end{proof}

Another useful way to state the preceding result is the following:
\begin{prop}
\label{Id-times-S-invariance}
Let $T$ and $S$ be two dynamical systems, with $S$ ergodic, and $\lambda\in J(T,S)$. If 
$\lambda$ is invariant by $\Id\times S$, then $\lambda=\mu\otimes\nu$.
\end{prop}

\subsubsection{Application: Furstenberg's multiple recurrence theorem for weakly mixing transformations}

As a nice application of Proposition~\ref{Id-times-S-invariance}, we give below an elegant proof of
Furstenberg's multiple recurrence theorem  in the (easy) case of weakly mixing transformations \cite{furst3}. 
This proof was given by V. Ryzhikov in \cite{ryzhi5}. Recall that $T$ is weakly mixing if $T\times T$ 
is ergodic, and that this implies 
\begin{itemize} 
\item $T\times S$ is ergodic for any ergodic system $S$ (see also subsection~\ref{subsection-product-ergodicity});
\item $T^k$ is weakly mixing for any $k\neq0$.
\end{itemize}

\begin{theo}
Let $T$ be a weakly mixing dynamical system. For any integer $k\ge1$, and any measurable subsets 
$A_0,\ldots,A_k$, we have 
\begin{equation}
\label{Furstenberg-WM}
\dfrac{1}{n}\sum_{j=1}^n \mu(A_0\cap T^{-j}A_1\cap T^{-2j}A_2\cap\cdots\cap T^{-kj}A_k)\ \tend{n}{\infty}\ \mu(A_0)\mu(A_1)\cdots\mu(A_k).
\end{equation}
\end{theo}
\begin{proof}
Consider for any $n\ge1$ the self-joining $\lambda_n$ of $T$ of order $k+1$ defined by
$$ \lambda_n(A_0\times A_1\times\cdots\times A_k)\ \egdef\ 
    \dfrac{1}{n}\sum_{j=1}^n \mu(A_0\cap T^{-j}A_1\cap T^{-2j}A_2\cap\cdots\cap T^{-kj}A_k). $$
Note that the result we want to prove is equivalent to the following translation in the language of joinings:
$$ \lambda_n\ \tend{n}{\infty}\ \mu\otimes\mu\otimes\cdots\otimes\mu. $$
By compacity, it is enough to verify that if $\lambda$ is any cluster point of the sequence $(\lambda_n)$,
then $\lambda$ is the product measure, which we shall prove by induction on $k$. Observe that such a cluster point 
is always $\Id\times T\times T^2\times\cdots\times T^{k}$-invariant. Thus Proposition~\ref{Id-times-S-invariance}
gives the result when $k=1$ by ergodicity of $T$. Now take $k>1$ and suppose the result is true for $k-1$. This
induction hypothesis ensures that the marginal of $\lambda$ on the last $k$ coordinates is the product measure.
Hence $\lambda$ can be seen as a joining of $\Id$ with $T\times T^2\times\cdots\times T^{k}$. Since 
 $T\times T^2\times\cdots\times T^{k}$ is ergodic (because $T$ is weakly mixing),
Proposition~\ref{Id-times-S-invariance} shows that the result is still true for~$k$. 
\end{proof}

\subsubsection{Disjointness and classification of systems}

As we have already mentioned, one of the reasons for the introduction of disjointness by Furstenberg
was the classification of classes of dynamical systems. And indeed there are many results in this direction:
For example, Furstenberg proved that a dynamical system has zero entropy if and only if it is disjoint from any Bernoulli dynamical system. It also follows from his work that $T$ is weakly mixing if and only if it is disjoint from any rotation on the circle. Some further explanations on how disjointness can characterize 
such classes of dynamical systems will be given in Section~\ref{stable-subsubsection} (see in particular Theorem~\ref{disjoint-classes-theorem}).  

\subsubsection{Disjointness and pointwise convergence of ergodic averages}
In~\cite{lesig2} was introduced the so-called \emph{weak disjointness} property of
dynamical systems:  
\begin{defi}
$(X,\A,\mu,T)$ and $(Y,\B,\nu,S)$ are said to be \emph{weakly disjoint} if, given any function $f$
in $L^{2}(\mu)$ and any function $g$ in $L^{2}(\nu)$, there exist a
set $A$ in $\A$ and a set $B$ in $\B$ such that
\begin{itemize}
           \item $\mu(A)=\nu(B)=1$
           \item for all $x\in A$ and for all $y\in B$, the sequence
           \begin{equation}\label{averages}
           \left(\frac1N\sum_{n=0}^{N-1}f\left(T^nx\right)\cdot
g\left(S^ny\right)
           \right)_{N>0}
          \end{equation}
           converges.
\end{itemize}
\end{defi}

(In fact it is sufficient to check this statement for only dense families of $f$
and $g$ in their respective $L^2$ spaces.) The main motivation for the introduction of 
this property was the study of non-conventional ergodic averages: Suppose that $S$ and 
$T$ act on the same probability space $X$, and consider for $f$ and $g$ in $L^2(\mu)$
the averages
$$\dfrac{1}{n}\sum_{k=0}^{n-1}f(T^kx)g(S^kx).$$
When do these averages converge $\mu$-a.e.?  If we assume that $T$ and $S$ are weakly 
disjoint, then almost-everywhere pointwise convergence holds. 

As we could guess from the name of the property, if $T$ and $S$ are disjoint,
they are weakly disjoint. The link between joinings and convergence of ergodic sums
like in \eqref{averages} comes from the following remark: We can always assume that $T$ and $S$ 
are continuous transformations of compact metric spaces $X$ and $Y$ respectively
(indeed, any measurable system is isomorphic to such a transformation on a compact metric space: see e.g. \cite{furst3}). Then the set of probability measures on $X\times Y$ equipped with the topology of weak convergence is metric compact. If moreover $T$ and $S$ are ergodic, we can easily find subsets $X_0\subset X$ and $Y_0\subset Y$ with $\mu(X_0)=\nu(Y_0)=1$, such that for all $(x,y)\in X_0\times Y_0$, any cluster point of the sequence
$$ \delta_n(x,y)\ =\ \dfrac{1}{n} \sum_{k=0}^{n-1} \delta_{(T^kx,S^ky)} $$
is automatically a joining of $T$ and $S$. When $T$ and $S$ are disjoint, there is therefore only
one possible cluster point to the sequence $\delta_n(x,y)$ which is $\mu\otimes \nu$. 
This ensures that for continuous functions $f$
and $g$, convergence of \eqref{averages} to the product of the integrals of $f$ and $g$ holds,
hence $T$ and $S$ are weakly disjoint.

It is not difficult to show that weak disjointness is really weaker than disjointness. Indeed,
there are many examples of dynamical systems which are self-weakly disjoint (weakly disjoint from themselves):
For example, irrational rotations and Chacon's transformation are self-weakly disjoint (see~\cite{lesig2}
for details and other examples). However,
a non trivial dynamical system is never disjoint from itself, as we will see in the next paragraph.

\subsection{Joinings supported on graphs}
\label{factor_joining_subsection}
Suppose that $S$ is a factor of $T$, with $\pi:\ X\to Y$ a homomorphism. Then $\pi$ gives rise
to the existence of a very special joining of $T$ and $S$, which we denote by $\Delta_\pi$, and which
is defined by
$$ \Delta_\pi(A\times B)\ \egdef\ \mu(A\cap\pi^{-1}B). $$

\begin{lemma}
\label{factor_joining_lemma}
The joining $\Delta_\pi$ has the following two properties:
\begin{itemize}
\item $\Delta_\pi(G)=1$, where $G:=\{(x,\pi(x)),\ x\in X\}$ is the graph of $\pi$;
\item $\{X,\emptyset\}\otimes\B\subset \A\otimes\{Y,\emptyset\}\ \mod \Delta_\pi$, where the notation
`$\,\C\subset \D\ \mod\lambda$' means that for each $C$ in the $\sigma$-algebra $\C$,
there exists a $D$ in the $\sigma$-algebra $\D$ such that $\lambda(C\Delta D)=0$.
\end{itemize}
\end{lemma}
\begin{proof}
Take a refining and generating sequence $(\P_n)$ of finite measurable partitions in $(Y,\B,\nu)$. 
(That is to say, $\P_{n+1}$ is finer than $\P_n$, and if $y$ and $y'$ are distinct points in $Y$, there
exists some $n$ for which $y$ and $y'$ do not belong to the same atom of $\P_n$; this ensures in particular 
that the $\sigma$-algebra generated by the partitions $\P_n$ is $\B$.)
Then $G$ can be written as
$$ G\ =\ \bigcap_n \bigcup_{B\in\P_n} \pi^{-1}(B)\times B. $$
But for any $n$, it is clear from the definition of $\Delta_\pi$ that
$$ \Delta_\pi\left( \bigcup_{B\in\P_n} \pi^{-1}(B)\times B\right) \ =\ 1, $$
which proves that $\Delta_\pi$ is supported on the graph of $\pi$. For the second property, 
it is easy to verify that $\Delta_\pi$ identifies any $X\times B$, where $B\in\B$, with
$\pi^{-1}B\times Y$.
\end{proof}

It is quite remarkable to see that the converse is true, \textit{i.e.} we can characterize
the fact that $S$ is a factor of $T$ by the existence of some special joining.

\begin{prop}
\label{factor_joining}
Suppose that there exists a joining $\lambda\in J(T,S)$ such that
$$ \{X,\emptyset\}\otimes\B\subset \A\otimes\{Y,\emptyset\}\ \mod \lambda. $$
Then $S$ is a factor of $T$ and we can find a homomorphism $\pi:\ X\to Y$ 
such that $\lambda$ is supported on the graph of $\pi$.
\end{prop}

\begin{proof}
Take again a refining and generating sequence $(\P_n)$ of finite measurable partitions in $(Y,\B,\nu)$. 
For each $n$, $\P_n$ can be identified via $\lambda$ with some $\overline{\P_n}$ in $X$,
where $(\overline{\P_n})$ is a refining, but not necessarily generating, sequence of partitions in $(X,\,\mu)$,
Now, for $\mu$-almost every $x\in X$, if we denote by $\overline{B_n}$ the atom of $\overline{\P_n}$ to which $x$ belongs and $B_n$ the atom of ${\P_n}$ identified via $\lambda$ with $\overline{B_n}$, then $\bigcap_n B_n$ is a singleton. Denote by $\pi(x)$ the only element of this set; we let the reader check that $x\mapsto \pi(x)$ satisfies the required properties.
\end{proof}

If we assume further that $T$ and $S$ are \emph{isomorphic} (in other words, the homomorphism $\pi$ from
$X$ to $Y$ is invertible and $\pi^{-1}$ is also a homomorphism of dynamical systems; in this case we will rather
use the symbol $\phi$ instead of $\pi$), the joining $\Delta_\phi$ now satisfies $\{X,\emptyset\}\otimes\B = \A\otimes\{Y,\emptyset\}\ \mod \Delta_\phi$, where the notation
`$\,\C = \D\ \mod\lambda$' means that we have both  $\,\C\subset \D\ \mod\lambda$ and $\,\D\subset \C\ \mod\lambda$.
The analog of Proposition~\ref{factor_joining} is also true, and we summarize all these results in the following
theorem.

\begin{theo}
\label{joinings-supported-on-graphs-thm}
Let $T$ and $S$ be two dynamical systems. $S$ is a factor of $T$ if and only if there exists some 
joining $\lambda\in J(T,S)$ such that $\{X,\emptyset\}\otimes\B\subset \A\otimes\{Y,\emptyset\}\ \mod \lambda$; in this case 
$\lambda$ is supported on the graph of some homomorphism. 
$T$ and $S$ are isomorphic if and only if there exists some 
joining $\lambda\in J(T,S)$ such that $\{X,\emptyset\}\otimes\B= \A\otimes\{Y,\emptyset\}\ \mod \lambda$; in this case 
$\lambda$ is supported on the graph of some isomorphism. 
\end{theo}

\subsubsection{Application: isomorphism of discrete-spectrum transformations}

There are very deep applications of the characterization of isomorphic systems by joinings. Krieger's theorem, 
stating that any dynamical system with entropy less than $\log_2(k)$ is isomorphic to the shift map
on $\{1,\ldots,k\}$ with an appropriate invariant measure, and Ornstein's theorem characterizing the dynamical 
systems which are isomorphic to some Bernoulli shift, can both be proved using this result (see \cite{rudol7}, Chapter~7). Here we illustrate the power of Theorem~\ref{joinings-supported-on-graphs-thm} by 
giving Lema\'nczyk's elegant proof of a well-known theorem due to Halmos and von Neumann. This proof is taken 
from \cite{thouv2}. Recall that a system $T$ is said to have \emph{discrete spectrum} if the subspace of $L^2$
spanned by the eigenvectors of $f\mapsto f\circ T$ is dense in $L^2$.

\begin{theo}
If $T$ is ergodic and has discrete spectrum, and if $S$ is spectrally isomorphic to $T$, then $T$ and $S$ are isomorphic.
\end{theo}

\begin{proof}
First note that $S$ is ergodic and also has discrete spectrum with the same eigenvalues as $T$, 
since it is spectrally isomorphic to $T$.
Let $\lambda$ be an ergodic joining of $T$ and $S$. For an eigenvalue $\alpha$ of $T$, let $f$ be
an eigenvector associated to $\alpha$ in the system $T$, and $g$ associated to $\alpha$ in the system $S$.
Then both functions $(x,y)\mapsto f(x)$ and $(x,y)\mapsto g(y)$ are eigenvectors in $(T\times S)_\lambda$.
Since $(T\times S)_\lambda$ is an ergodic system, $\alpha$ is a simple eigenvalue in  $(T\times S)_\lambda$, 
hence  $(x,y)\mapsto g(y)$ and $(x,y)\mapsto f(x)$ are proportional. This implies that $(x,y)\mapsto f(x)$
belongs to $L^2(\{X,\emptyset\}\otimes\B,\lambda)$; but since eigenvectors are dense in $L^2(\A,\mu)$, this in turn 
implies $L^2(\A\otimes\{Y,\emptyset\},\lambda)\subset L^2(\{X,\emptyset\}\otimes\B,\lambda)$. For the same reason, the inverse inclusion
is also true, which proves that $\{X,\emptyset\}\otimes\B= \A\otimes\{Y,\emptyset\}\ \mod \lambda$, and the two systems are isomorphic.
\end{proof}

\subsubsection{Self-joinings, commutant and minimal self-joinings}

When studying self-joinings of a single dynamical system $T$, we are obviously in a situation
where the two systems that we want to join together are isomorphic. Now, thinking of what 
an isomorphism between $T$ and itself is, we see that this has to be an automorphism of the Lebesgue space
$(X,\A,\mu)$, which commutes with $T$.
We call \emph{commutant of $T$}, and denote by $C(T)$, the set
$$ C(T)\ \egdef\ \{S\in\Aut(X,\A,\mu): S\circ T = T\circ S\}. $$
For each $S\in C(T)$, let us denote by $\Delta_S$ the self-joining of $T$ defined by
$$ \Delta_S(A\times B)\ \egdef\ \mu(A\cap S^{-1}B). $$
This leads to an identification of  $C(T)$ with some subset of $J_2(T)$, namely the subset of 
all $\lambda\in J_2(T)$ such that 
$$ \A\otimes \{X,\emptyset\}\ =\ \{X,\emptyset\}\otimes\A\ \mod \lambda. $$
Observe that for each $T$ the commutant of $T$ contains at least the powers $T^n$, $n\in\ZZ$. Besides,
note that each joined system $(T\times T)_{\Delta_S}$ for $S\in C(T)$ is isomorphic to $T$ (an isomorphism
is given by $x\longmapsto(x,Sx)$). In particular,  $(T\times T)_{\Delta_S}$ is ergodic as soon as $T$ is
ergodic. This gives for each ergodic $T$ a family of ``obvious'' ergodic self-joinings: the $\Delta_{T^n}$,
$n\in\ZZ$. In the case where $T$ is weakly mixing, the product measure $\mu\otimes\mu$ is also considered 
as an obvious ergodic self-joining. In 1979, Rudolph introduced in \cite{rudol10} the important notion of ``minimal self 
joinings'', which says that $T$ has no other ergodic self-joinings of order 2 than the obvious ones. 
The  notion of ``$k$-fold minimal self-joinings'' refers to joinings of order $k$:
\begin{defi}
For $k\ge 2$, we say that an ergodic $T$ has \emph{$k$-fold minimal self-joinings} if for
each ergodic joining $\lambda$ of $k$ copies of $T$, we can partition the set $\{1,\ldots,k\}$ of coordinates
into two subsets $J_1,\ldots,J_\ell$ such that
\begin{enumerate}
\item \label{first-condition}
for $j_1$ and $j_2$ belonging to the same $J_i$, the marginal of $\lambda$ on the coordinates $j_1$ and $j_2$ is some  $\Delta_{T^n}$;
\item for $j_1\in J_1,\ldots,j_\ell\in J_\ell$, the coordinates $j_1\ldots,j_\ell$ are independent.
\end{enumerate}
We say that $T$ has minimal self-joinings (MSJ) if $T$ has $k$-fold minimal self-joinings for every $k\ge2$.
\end{defi}

In fact, it has been proved by Glasner, Host and Rudolph \cite{glasn2} that for a weakly mixing $T$, 
3-fold minimal self-joinings implies $k$-fold minimal self-joinings for all $k$. (See also a proof in \cite{thouv2}.)
A very important question which remains open in the area is the following:
\begin{question}
 Does there exist a transformation $T$ which has 2-fold but not 3-fold minimal self-joinings?
\end{question}
 (This question is closely related to the problem of
2-fold and 3-fold mixing: See Section~\ref{subsection-2-3-mixing}.)

There are many examples of dynamical systems with the MSJ property, the simplest of which is 
probably the well known Chacon's transformation (see~\cite{delju5}). One immediate consequence of this property
is that $C(T)$ is reduced to the powers of $T$. This in turn implies that the transformation $T$ has no square root
(\textit{i.e. there is no automorphism $S$ such that $S\circ S=T$}). It can be noticed that the first example
of a transformation with no square root, given by Ornstein \cite{ornst1} in 1970, belongs to the class of
\emph{mixing rank-one\footnote{We shall not define in this text what a rank-one system is; the reader unfamiliar
with this notion can find a good introduction in \cite{nadka1}} systems}, 
which all turned out to have the MSJ property, as King proved in 1988~\cite{king2}.
Another joining proof of this result has been given by Ryzhikov~\cite{ryzhi1}.

\medskip

If, in condition~(\ref{first-condition}) of the definition of MSJ, we replace $\Delta_{T^n}$ 
by $\Delta_S$, for some $S\in C(T)$, we get the weaker, but important, notion of \emph{simplicity}.
In the case of joinings of order~2, this property was introduced by Veech~\cite{veech1}, who proved a
useful result on the structure of factor $\sigma$-algebras of 2-simple systems. The general definition of
simplicity was later given by Del Junco and Rudolph~\cite{delju3}. We also refer to Thouvenot~\cite{thouv2}
for a presentation of properties of simple systems.
 
\medskip

As a nice application of joinings to the study of $C(T)$, we should also mention here Ryzhikov's proof 
of King's so-called \emph{Weak Closure Theorem}, stating that if $T$ is rank one, any $S$ in $C(T)$ is 
the limit of some sequence $T^{n_k}$ of the powers of $T$ \cite{king1,ryzhi1}. 


\section{Joinings and factors}\label{joinings-factors-section}

In this section we discuss the links and the differences between the two properties proposed in the introduction
to define ``$T$ and $S$ are relatively prime''. The starting point of the theory is the existence of a special joining 
arising when $T$ and $S$ have a common factor. 

\subsection{Relatively independent joining above a common factor}

We first define this special joining in the particular case 
of two systems which obviously have a common factor, since they both arise as a joining of a fixed system $R$ with other systems. 

Let $(X,\A,\mu,T)$, $(Y,\B,\nu,S)$ and $(Z,\C,\rho,R)$ be three dynamical systems. Suppose that $\lambda_T$ is a joining of
$T$ and $R$, and $\lambda_S$ is a joining of $S$ and $R$. Then we can put all these systems together so that the marginal
on $(X\times Z)$ is $\lambda_T$ and the marginal on $(Y\times Z)$ is $\lambda_S$: First pick $z$ according to the probability 
law $\rho$, then pick $x$ and $y$ according to their conditionnal law knowing $z$ in the respective joinings $\lambda_T$ and 
$\lambda_S$, but independently of each other. More precisely, consider the 3-fold joining of $T$, $S$
and $R$ denoted by $\lambda_T\otimes_R\lambda_S$ defined by setting, for all $A\in \A$, $B\in\B$ and $C\in\C$
\begin{equation}
\lambda_T\otimes_R\lambda_S(A\times B\times C)\ \egdef\ \int_C \EE_{\lambda_T}[\ind{x\in A}|z] \, \EE_{\lambda_S}[\ind{y\in B}|z] \, d\rho(z).
\end{equation}
We get in this way a joining of the two systems $(T\times R)_{\lambda_T}$ and $(S\times R)_{\lambda_S}$ identifying 
the coordinate $z$ on both systems. This joining is called the \emph{relatively independent joining of $(T\times R)_{\lambda_T}$
and $(S\times R)_{\lambda_S}$ above their common factor $R$}.

Now turn to the more general situation of two systems $T$ and $S$ sharing a common factor $R$. 
This means that we have homomorphisms of dynamical systems $\pi_T:\ X\to Z$ and 
$\pi_S:\ Y\to Z$, which give us special joinings $\lambda_{\pi_T}\in J(T,R)$ and $\lambda_{\pi_S}\in J(S,R)$ (see
section~\ref{factor_joining_subsection}). We can then consider the relatively independent joining $\lambda_{\pi_T}\otimes_R
\lambda_{\pi_S}$, which is a 3-fold joining of $T$, $S$ and $R$. The marginal of $\lambda_{\pi_T}\otimes_R\lambda_{\pi_S}$
on $X\times Y$ is now a joining of $T$ and $S$, called the \emph{relatively independent joining of $T$
and $S$ above their common factor $R$}, and denoted by $\mu\otimes_R\nu$. This joining 
has the strong property of identifying the projections of both systems on the common factor.

In summary, the $\sigma$-algebra $\C$ can be viewed as a sub-$\sigma$-algebra of both $\A$ and $\B$, any $\C$-measurable function can be viewed as a function defined on $Z$, and we have the following integration formula:
$$
\int_{X\times Y} f(x)\,g(y)\, d(\mu\otimes_R\nu)(x,y)\ =\ \int_Z\EE[f(x)|\C](z)\,\EE[g(y)|\C](z)\,d\rho(z).
$$
\begin{figure}
\begin{picture}(0,0)%
\includegraphics{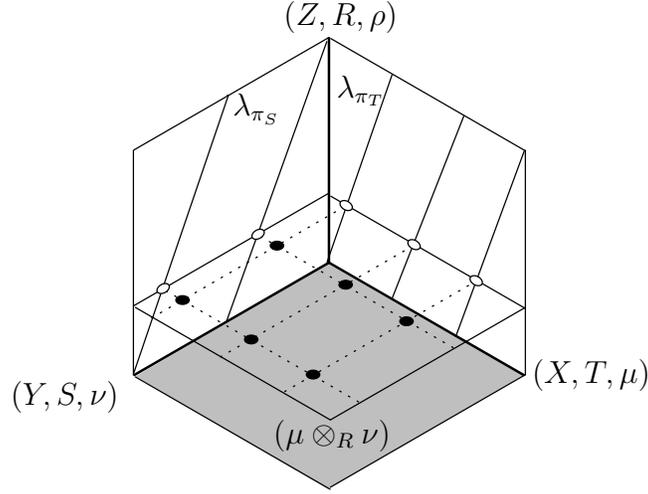}%
\end{picture}%
\setlength{\unitlength}{4144sp}%
\begingroup\makeatletter\ifx\SetFigFont\undefined%
\gdef\SetFigFont#1#2#3#4#5{%
  \reset@font\fontsize{#1}{#2pt}%
  \fontfamily{#3}\fontseries{#4}\fontshape{#5}%
  \selectfont}%
\fi\endgroup%
\begin{picture}(5402,3152)(225,-2627)
\put(1021,-2041){\makebox(0,0)[lb]{\smash{{\SetFigFont{12}{14.4}{\rmdefault}{\mddefault}{\updefault}{\color[rgb]{0,0,0}$(Y,S,\nu)$}%
}}}}
\put(2595,-2268){\makebox(0,0)[lb]{\smash{{\SetFigFont{12}{14.4}{\rmdefault}{\mddefault}{\updefault}{\color[rgb]{0,0,0}$(\mu\otimes_R\nu)$}%
}}}}
\put(4141,-1906){\makebox(0,0)[lb]{\smash{{\SetFigFont{12}{14.4}{\rmdefault}{\mddefault}{\updefault}{\color[rgb]{0,0,0}$(X,T,\mu)$}%
}}}}
\put(2971,-196){\makebox(0,0)[lb]{\smash{{\SetFigFont{12}{14.4}{\rmdefault}{\mddefault}{\updefault}{\color[rgb]{0,0,0}$\lambda_{\pi_T}$}%
}}}}
\put(2351,-308){\makebox(0,0)[lb]{\smash{{\SetFigFont{12}{14.4}{\rmdefault}{\mddefault}{\updefault}{\color[rgb]{0,0,0}$\lambda_{\pi_S}$}%
}}}}
\put(2656,238){\makebox(0,0)[lb]{\smash{{\SetFigFont{12}{14.4}{\rmdefault}{\mddefault}{\updefault}{\color[rgb]{0,0,0}$(Z,R,\rho)$}%
}}}}
\end{picture}%
\caption{The relatively independent joining of $T$ and $S$ above their common factor $R$: Given
a (common) fixed projection on the factor, the joining is the direct product of the two conditional measures.}
\end{figure} 

\begin{prop}
\label{RIJ-prop}
Let $\mu{\otimes}_{R}\nu$ be the relatively independent joining of $T$ and $S$  above their common
factor $R$. Then 
$$ \pi_T(x)\ =\ \pi_S(y)\quad \mu{\otimes}_{R}\nu\text{-a.e.}$$
\end{prop}
\begin{proof}
By construction, $\mu{\otimes}_{R}\nu$ is the marginal on $X\times Y$ of the 3-fold joining 
$\lambda_{\pi_T}\otimes_R\lambda_{\pi_S}$ of $T$, $S$ and $R$
whose marginal on $X\times Z$ (respectively $Y\times Z$) is $\lambda_{\pi_T}$ (respectively $\lambda_{\pi_S}$). 
Recall that, from Lemma~\ref{factor_joining_lemma},  we have $z=\pi_T(x)$
$\lambda_{\pi_T}$-a.e. and  $z=\pi_S(y)$ $\lambda_{\pi_S}$-a.e. Therefore, $z=\pi_T(x)=\pi_S(y)$
$\lambda_{\pi_T}\otimes_R\lambda_{\pi_S}$-a.e., hence
$\pi_T(x)=\pi_S(y)$ $\mu{\otimes}_{R}\nu$-a.e.
\end{proof}

Note that what we did here with two systems could have been done with any finite or countable family
of systems sharing a common factor.

\subsection{Self-joinings and factors}

Any factor $S$ of $T$ can be considered as a common factor of two copies of $T$, thus 
gives rise to a special self-joining of $T$: The relatively independent self-joining of $T$ above its factor $S$.
As in ordinary arithmetic of integer numbers, every dynamical system which is not reduced to one point has at least two factors:
$T$ itself is of course a factor of $T$, and the system reduced to one point is also a factor of $T$. 
Note that the relatively independent joining of $T$ above $T$ itself
is the self-joining $\Delta_{\text{Id}}$ of $T$ supported on the graph of $\Id$, 
and that the relatively independent joining of $T$ above the 
factor reduced to one point is the product measure $\mu\otimes\mu$. Now, we want to show that any other factor of $T$
is incompatible with the 2-fold MSJ property.

\begin{prop}
Let $(X,\A,\mu,T)$ be an ergodic dynamical system where $X$ is not finite. 
If $T$ has another non-trivial factor than $T$ itself, then $T$ does not have 
the 2-fold MSJ property.
\end{prop}
\begin{proof}
Let $S$ be a non-trivial factor of $T$, and consider $\mu\otimes_S\mu$ the relatively independent self-joining of
$T$ above $S$ constructed from the homomorphism $\pi:\ X\to Y$. 
The system defined by this joining may not be ergodic, but if we assume further that $T$ has the 
2-fold MSJ property, the ergodic decomposition of $\mu\otimes_S\mu$ has the form
$$ \mu\otimes_S\mu\ =\ \theta \mu\otimes \mu + (1-\theta)\sum_{n\in\ZZ}p_n\Delta_{T^n}, $$
where $0\le\theta\le 1$, $p_n\ge 0$ and $\sum_{n\in\ZZ}p_n=1$. Denote by $x$ and $x'$ the two coordinates 
in the Cartesian square of $X$, and remember that, by Proposition~\ref{RIJ-prop}, we have $\pi(x)=\pi(x')$ 
$\mu\otimes_S\mu$-a.e. Since $S$ is not the factor reduced to one point, this would not be possible if 
$\mu\otimes\mu$ appeared in the ergodic decomposition, hence $\theta=0$. Moreover, if the ergodic decomposition
of $\mu\otimes_S\mu$ was reduced to $\Delta_{\text{Id}}$ , $S$ would be $T$ itself. Hence 
there is some $n\neq 0$ such that $p_n>0$. In other words, there exists some $n\neq 0$ such that, under $\mu\otimes_S\mu$,
\begin{equation}
\mu\otimes_S\mu(x'=T^nx)\ >\ 0.
\end{equation}
This implies that, in the factor $S$, $y=S^ny$ with positive probability.
But $S$ is ergodic (because $T$ is), hence the space $Y$ on which $S$ acts is finite. In turn, this implies
that there exists at least one $y\in Y$ such that 
\begin{equation}
\label{atomic}
\mu\otimes_S\mu\Bigl((x'=T^nx)\cap(\pi(x)=y)\Bigr)\ >\ 0.
\end{equation}
Now, observe that $\mu\otimes_S\mu$ conditioned on $(\pi(x)=y)$ is the direct product of $\mu$ conditioned on $(\pi(x)=y)$
with itself. Thus, \eqref{atomic} implies that $\mu$ conditioned on $(\pi(x)=y)$ has at least one atom. But $\mu(\pi(x)=y)>0$,
hence $\mu$ itself has an atom, hence $T$ is periodic.
\end{proof}

\begin{corollary}
If $T$ is ergodic aperiodic with the 2-fold MSJ property, then any non constant stationary process
living in the system $T$ generates the whole $\sigma$-algebra.
\end{corollary}

\subsection{Disjointness and lack of common factor}
\label{dlcf-subsection}
We now turn back to the case of two dynamical systems $T$ and $S$ sharing a common factor $R$. 
As we have already noticed in the case of self-joinings, Proposition~\ref{RIJ-prop} ensures that if
$R$ is not reduced to one point, the relatively independent joining of $T$ and $S$ above $R$ is not 
the product measure $\mu\otimes\nu$, hence $T$ and $S$ are not disjoint. This means that Property~\ref{disjointness-property}
(the disjointness property) implies Property~\ref{no-common-factor-property} (the lack of common factor). It had been asked by Furstenberg \cite{furst1} whether these two properties were, in fact equivalent. The discovery by Rudolph of systems with the MSJ property allowed him to answer negatively to this question~\cite{rudol10}.

\begin{prop}
There exist two non disjoint systems $S$ and $T$ without any other common factor than the trivial one.
\end{prop}
\begin{proof}
Let $\xi=(\xi_k)$ be a stationary process generating a system $T$ with the MSJ property, and let $(\xi'_k)$ be an independent
copy of $(\xi_k)$. We now define another stationary process $\zeta$ by setting 
$$ \zeta_k\ \egdef\ \{\xi_k,\xi'_k\}. $$
(This means that $\zeta_k$ tells us which are the two values taken by $\xi_k$ and $\xi'_k$, but neither which one is 
taken by $\xi_k$ nor which one is taken by $\xi'_k$.) Denote by $S$ the system generated by $\zeta$: $S$ is a factor
of the Cartesian product $T\times T$. Then $T$ and $S$
are certainly not disjoint, since the two processes $\xi$ and $\zeta$ are clearly not independent of each other.
However, suppose that $S$ and $T$ share a common factor which is not the system reduced to one point. Then, 
since $T$ has the MSJ property, this factor has to be $T$ itself. Therefore, we can find a third copy $\xi''$ of $\xi$, living in
the  Cartesian square $T\times T$, and 
which is measurable with respect to the factor $\sigma$-algebra generated by $\zeta$. By the 3-fold 
MSJ property of $T$, either $\xi''$ is equal to one of the processes $\xi$ or $\xi'$, possibly shifted by some constant integer
$k$, or $\xi''$ is independent of $(\xi,\xi')$. The former case is impossible since $\zeta$, hence $\xi''$, is invariant
by the ``flip'' $(x,x')\mapsto(x',x)$, which is clearly not true for any shifted copy of $\xi$ or $\xi'$. The latter case is also
impossible since $(\xi)$ and $(\xi')$ generate the whole $\sigma$-algebra in $T\times T$.
\end{proof}

\subsubsection{The fundamental lemma of non-disjointness}

There exists however an important result which tells us how to derive some information on factors
from the non-disjointness of two systems. This theorem appeared for the first time in year~2000, in two publications \cite{glasn3,leman6}.

\begin{theo}
\label{non-disjointness}If $T$ and $S$ are not disjoint, then $S$ has a non trivial common factor
with some joining of a countable family of copies of $T$.
\end{theo}
\begin{proof}
Start with a joining $\lambda$ of $(X,\A,\mu,T)$ and $(Y,\B,\nu,S)$: Then $S$ is a factor of $(T\times S)_\lambda$. Consider a 
countable family of copies of $(T\times S)_\lambda$, and denote by $\lambda_\infty$ their relatively independent 
joining over their common factor $S$. Since this joining identifies all the projections on $Y$, it can be viewed as a
probability measure on $Y\times X^\NN$, and is characterized by
$$\lambda_\infty(B\times A_0\times\cdots\times A_k\times X\times X\cdots)\ =\ 
\int_B \EE_\lambda[\ind{A_0}|y]\cdots\EE_\lambda[\ind{A_k}|y]\,d\nu(y). $$
Observe that this probability measure is invariant under the shift on each $y$-fiber:
$$ (y,x_0,x_1,\ldots)\ \mathop{\longmapsto}\ (y,x_1,x_2\ldots). $$
Moreover, $\lambda_\infty$ conditioned on such a fiber is a product measure: The infinite direct
product of $\lambda$ conditioned on $y$. A relative version of Kolmogorov 0-1 law 
(precisely stated in \cite{leman6}, Lemma~9) then tells us that the 
$\sigma$-algebra of shift-invariant events coincides modulo $\lambda_\infty$ with the
$\sigma$-algebra of $y$-measurable events. 

Now, let $f$ be a bounded measurable function defined on $X$, which we also consider as a function
on $X\times Y$ and on $Y\times X^\NN$ by setting $f(x,y)\egdef f(x)$, and $f(y,x_0,x_1,\ldots)\egdef f(x_0)$.
Applying Birkhoff ergodic theorem to $f$ in the system defined by the shift and $\lambda_\infty$ gives
\begin{equation}
\label{birkhoff}\lim_{n\to\infty} \dfrac{1}{n}\sum_{k=0}^{n-1} f(x_k)\ =\ \EE_{\lambda_\infty}[f|y]\ =\ \EE_{\lambda}[f|y] 
\quad \lambda_\infty\text{-a.e.}
\end{equation}
Since $T$ and $S$ are not disjoint, we could have chosen $\lambda\neq\mu\otimes\nu$, and then $f$ such that
$\EE_{\lambda}[f|y]$ is not constant modulo $\lambda$. Then this conditional expectation generates in $S$
a non-trivial factor, which by \eqref{birkhoff} can be identified to some factor of the joining of
a countable family of copies of $T$ defined by the marginal of $\lambda_\infty$ on $X^\NN$.
\end{proof}

\subsubsection{Application: disjointness of classes of dynamical systems}
\label{stable-subsubsection}
We get from Theorem~\ref{non-disjointness} some important disjointness 
results.
Note that some properties of dynamical systems are stable under the operations of taking
joinings and factors. We will call these properties \emph{stable properties}.
This is \textit{e.g.} the case of the zero-entropy property: We know that
a factor of a zero-entropy system still has zero entropy, and that any joining of zero-entropy
systems also has zero entropy. Take now any $K$-system $S$: We know that any non-trivial factor 
of $S$ has positive entropy, hence $S$ cannot have a non-trivial common factor with a joining of 
copies of a zero-entropy $T$.

The same argument also applies to the disjointness of discrete-spectrum systems with weakly mixing
systems, since discrete spectrum is a stable property, and weakly mixing systems are 
characterized by the fact that they do not have any discrete-spectrum factor.
We thus have the following theorem:

\begin{theo}
\label{disjoint-classes-theorem}
Any zero-entropy system is disjoint from any $K$-system. Any discrete-spectrum system is disjoint
from any weakly mixing system.
\end{theo}

\subsection{Joinings and $T$-factors}
\label{T-factor-subsection}
The result of Theorem~\ref{non-disjointness} leads to the introduction of a special class
of factors when some dynamical system $T$ is given: For any other dynamical system $S$,  
call \emph{$T$-factor} of $S$ any common factor of $S$ with
a joining of countably many copies of $T$. We also call \emph{$T$-factor $\sigma$-algebra}
any factor $\sigma$-algebra of $S$ associated to a $T$-factor of $S$.
Another way to state Theorem~\ref{non-disjointness}
is then the following: If $S$ and $T$ are not disjoint, then $S$ has a non-trivial
$T$-factor. In fact, the proof of this theorem gives a more precise result:
For any joining $\lambda$ of $S$ and $T$, for any bounded measurable
function $f$ on $X$, the factor $\sigma$-algebra of $S$ generated by the function
$\EE_{\lambda}[f(x) | y]$ is a $T$-factor $\sigma$-algebra of
$S$. 

With the notion of $T$-factor, we can extend Theorem~\ref{non-disjointness} in the following way,
showing the existence of a special $T$-factor $\sigma$-algebra of $S$ concentrating anything in $S$
which could lead to a non trivial joining between $T$ and $S$. 

\begin{theo}
\label{T-fact-max}
Given two dynamical systems $(X,\A,\mu,T)$ and $(Y,\B,\nu,S)$, there always exists
a maximum $T$-factor $\sigma$-algebra of $S$, denoted by $\F_T$.

Under any joining $\lambda$ of $T$
and $S$, the $\sigma$-algebras $\A\otimes\{\emptyset,Y\}$ and
$\{\emptyset,X\}\otimes\B$ are independent conditionally to the
$\sigma$-algebra $\{\emptyset,X\}\otimes\F_T$.
\end{theo}

The proof of the theorem is based on the two following lemmas.
\begin{lemma}\label{Tfm1}
Let $(B_i)_{i\in I}$ be a countable family of events in $\B$, such that
for all $i$, $B_i$ belongs to some $T$-factor $\sigma$-algebra $\F_i$
of $S$. Then there exists a $T$-factor $\sigma$-algebra of $S$ containing
all the $B_i$'s.
\end{lemma}
\begin{proof}
For each $i\in I$, we have a joining $\lambda_i$ of $S$ with a countable
family $(T_{i,n})_{n\in\NN}$ of copies of $T$, such that
$\F_i\subset\bigotimes_{n\in\NN}\A_{i,n} \mod\lambda_i$. Let us denote by
$R_i$ the dynamical system defined by $\lambda_i$, and by $\lambda$ the
relatively independent joining of all the $R_i$, $i\in I$, over their
common factor $S$. We can see $\lambda$ as a joining of $S$ with the
countable family $(T_{i,n})_{(i,n)\in I\times\NN}$ and, for each $i$ we
have
$$ \F_i\ \subset\ \bigotimes_{(i,n)\in I\times\NN}\A_{i,n}\quad\mod
\lambda. $$
We conclude that the factor $\sigma$-algebra of $S$ generated by all the $\F_i$'s is a
$T$-factor $\sigma$-algebra, which certainly contains all the $B_i$'s.
\end{proof}

\begin{lemma}\label{Tfm2}
Let $\F$ be a factor $\sigma$-algebra of $S$. If there exists a joining $\lambda$ of $T$
and $S$ under which the $\sigma$-algebras $\A\otimes\{\emptyset,Y\}$ and
$\{\emptyset,X\}\otimes\B$ are not independent conditionally to
$\{\emptyset,X\}\otimes\F$, then there exists a $T$-factor $\sigma$-algebra $\F'$ of $S$,
not contained in $\F$.
\end{lemma}
\begin{proof}
The hypothesis of the lemma implies the existence of a bounded measurable
function $f$ on $X$ such that, on a set of positive $\nu$-measure,
$$ \EE_{\lambda}[f(x) | \B]\ \not=\ \EE_{\lambda}[f(x) | \F]. $$
The factor $\sigma$-algebra $\F'$ of $S$ generated by the function $\EE_{\lambda}[f(x) | \B]$
is not contained in $\F$; but we saw in the proof of Theorem
\ref{non-disjointness} that $\F'$ is a $T$-factor $\sigma$-algebra.
\end{proof}

\begin{proof}[Proof of Theorem \ref{T-fact-max}]
In order to prove the existence of a maximum $T$-factor $\sigma$-algebra, we define
$$ \F_T\ \egdef\
          \{B\in \B\ :\ \mbox{$B$ belongs to a $T$-factor $\sigma$-algebra of $S$}\}\;, $$
and we claim that it is a $T$-factor $\sigma$-algebra. Since $(Y,\B,\nu)$ is a Lebesgue
space, the $\sigma$-algebra $\B$ equipped with the metric
$d(B,C):= \nu(B\Delta C)$ is separable (as usual, we identify subsets
$B$ and $C$ of $Y$ when $\nu(B\Delta C)=0$). There exists a countable
family $(B_i)_{i\in I}$ dense in $\F_T$, and, thanks to Lemma \ref{Tfm1},
there exists a $T$-factor $\sigma$-algebra $\F$ containing all the $B_i$'s. By density, we
have $\F_T\subset\F$ but, since $\F_T$ contains all the $T$-factor $\sigma$-algebras, we
have $\F_T=\F$.
This proves the first assertion of Theorem \ref{T-fact-max}. The second one
is just the application of Lemma \ref{Tfm2} to $\F_T$.
\end{proof}

The notion of $T$-factor was introduced in \cite{lesig2} in order to study
some aspects of the weak disjointness of dynamical systems. 
It was used there to prove that if $T$ satisfies the property of being 
\emph{self-weakly disjoint to all order $k\ge2$} (which is the natural generalization
of self weak-disjointness to the case of $k$ copies of $T$), then $T$ is
weakly disjoint from \emph{any other dynamical system.} The main open question
in this area is whether the self weak disjointness of $T$ alone is sufficient to
ensure this universal weak disjointness property.

\subsubsection*{A characterization of disjointness?}

We know from Theorem~\ref{non-disjointness} that if $T$ and $S$ are not disjoint, then
the maximum $T$-factor of $S$ is not trivial. One may ask whether the converse is true, but
in fact it is not difficult to find a counterexample: Indeed, take $(\sigma_n)_{n\in\ZZ}$  any stationary
process taking its values in $\{0,1\}$, and let $(\tau_n)_{n\in\ZZ}$ be an i.i.d. process
with $\tau_n$ taking each value 0 or 1 with probability $1/2$, independent of $(\sigma_n)$. Then, setting
$$ \tilde{\tau_n}\ \egdef\ \tau_n + \sigma_n \quad\mod 2, $$
we get another i.i.d. process $(\tilde{\tau}_n)$ with the same distribution as $(\tau_n)$.
The $(\sigma_n)$ process is clearly measurable with respect to the $\sigma$-algebra generated
by $(\tau_n)$ and $(\tilde{\tau}_n)$, which means that the dynamical system $S$ generated by
$(\sigma_n)$ is a factor of a self-joining of the Bernoulli shift generated by $(\tau_n)$. Therefore,
denoting by $T$ this Bernoulli shift, we get that any dynamical system $S$ generated by a two-valued
stationary process is a $T$-factor. But we could have taken for $(\sigma_n)$ any zero-entropy process,
which is disjoint from any Bernoulli shift (Theorem~\ref{disjoint-classes-theorem}).

\medskip

However, we can notice that Theorem~\ref{non-disjointness} gives a stronger conclusion when $T$
and $S$ are not disjoint. Indeed, this hypothesis being symmetric in $S$ and $T$, we get in that case 
that both the maximum $T$-factor of $S$ and the maximum $S$-factor of $T$ are not trivial. In the counterexample
given above,  $S$ having zero entropy, the maximum $S$-factor of a Bernoulli shift is trivial. So we ask the 
following question: 

\begin{question}
Can we find two disjoint systems $S$ and $T$ such that both the maximum $T$-factor of $S$ 
and the maximum $S$-factor of $T$ are not trivial?
\end{question}

\medskip

\subsubsection*{A remark on stable properties}

Analyzing the proof of Theorem~\ref{T-fact-max}, we can observe that the existence of
the maximum $T$-factor $\sigma$-algebra is in fact a corollary of the following: Being
a $T$ factor is a \emph{stable property} (in the sense given in Section~\ref{stable-subsubsection}).
And indeed, instead of the notion of $T$-factor we could take any stable property and we would
get by the same argument the existence of a maximum factor $\sigma$-algebra satisfying this 
stable property. Some classical factor
$\sigma$-algebras are particular cases of this type: If we start from the stable property of having 
zero entropy, we get the \emph{Pinsker $\sigma$-algebra} of the system, and if we start from the 
discrete-spectrum property, we get the \emph{Kronecker factor $\sigma$-algebra}.

\section{Self-Joinings and mixing}

We now turn to the relationships between the study of joinings and some
(weak and strong) mixing properties. We begin by a beautiful proof by
Ryzhikov of a well-known result concerning weak mixing.

\subsection{Weak mixing and ergodicity of products}\label{subsection-product-ergodicity}

There are many different (but equivalent!) ways to define the property of weak mixing.
A very concise one is the following: We say that $T$ is weakly mixing if the Cartesian
square $T\times T$ is ergodic. It is a standard result in ergodic theory that this
property implies the seemingly stronger one: 

\begin{theo}
\label{weak-mixing-thm}
Suppose that $T\times T$ is ergodic. Then for any ergodic $S$, $T\times S$ is ergodic.
\end{theo}

The classical proof of this implication makes use of the spectral theory of the action
of the unitary operators $U_T:\ f\mapsto f\circ T$ and 
$U_S:\ g\mapsto g\circ S$ on the subspaces of $L^2(\mu)$ and $L^2(\nu)$ of functions
with zero mean. It can be shown that if $T$ and $S$ are ergodic, 
the Cartesian square $T\times S$ is ergodic if and only if $U_T$ and 
$U_S$  do not have any common eigenvalue.

We propose here an alternative approach by joinings to prove Theorem~\ref{weak-mixing-thm}, 
due to Ryzhikov \cite{ryzhi9}. 
Let us start with a simple lemma on relatively independent joinings.
\begin{lemma}
\label{symmetrization-Lemma}
Let $\lambda$ be a joining of $T$ and $S$, and let $\lambda\otimes_S\lambda$ be the relatively independent
self-joining of $(T\times S)_\lambda$ above its factor $S$. If the marginal of $\lambda\otimes_S\lambda$
on $X\times X$ is the product measure $\mu\otimes\mu$, then $\lambda$ is the product measure $\mu\otimes\nu$. 
\end{lemma}
\begin{proof}
For all $A\in\A$, we have by definition of $\lambda\otimes_S\lambda$
$$ \lambda\otimes_S\lambda(A\times Y\times A)\ =\ \int_Y \left(\EE_\lambda[\ind{A}|y]\right)^2 \, d\nu(y). $$
Thus, if the marginal on $X\times X$ is the product measure, we get that for all $A\in\A$, 
$$ \left(\mu(A)\right)^2\ =\ \int_Y \left(\EE_\lambda[\ind{A}|y]\right)^2 \, d\nu(y)\ =\ 
   \left(\int_Y\EE_\lambda[\ind{A}|y] \, d\nu(y)\right)^2, $$
which is possible only if $\EE_\lambda[\ind{A}|y]$ is a constant $\nu$-a.e., hence 
$\lambda=\mu\otimes\nu$.
\end{proof}

\begin{proof}[Proof of Theorem~\ref{weak-mixing-thm}]
Assuming that $\mu\otimes\mu$ is an ergodic self-joining of $T$, we are 
going to show that for all ergodic $S$, $\mu\otimes\nu\in J_e(T,S)$, by
checking that it is an extremal point in $J(T,S)$, and applying 
Proposition~\ref{characterization-of-ergodic-joinings}.

Suppose that for some $0<p<1$ and some joinings $\lambda_1$ and $\lambda_2$ of $T$
and $S$, we have
$$ \mu\otimes\nu\ =\ p\lambda_1 + (1-p)\lambda_2. $$
Note that the relatively independent self-joining of $\mu\otimes\nu$ above the
factor $S$ is nothing but $\mu\otimes\nu\otimes\mu$, which can be decomposed as follows:
$$\mu\otimes\nu\otimes\mu\ =\ p^2 \lambda_1\otimes_S\lambda_1
                            + (1-p)^2 \lambda_2\otimes_S\lambda_2
			    +p(1-p) \lambda_1\otimes_S\lambda_2
			    +p(1-p) \lambda_2\otimes_S\lambda_1.
$$
Since $\mu\otimes\mu$ is ergodic, each term of the RHS must have $\mu\otimes\mu$ as
marginal on $X\times X$. But this implies that $\lambda_1$ and $\lambda_2$ satisfy the
hypothesis of Lemma~\ref{symmetrization-Lemma}, hence both $\lambda_1$ ad $\lambda_2$
are the product measure $\mu\otimes\nu$.
\end{proof}

\subsection{Characterization of mixing by joinings}

Strong mixing is also very easily characterized in term of joinings: By definition, the dynamical
system $T$ is \emph{strongly mixing} (or simply \emph{mixing}) if $\forall A,B\in\A$,
$$\mu(A\cap T^{-n}B)\tend{n}{\infty} \mu(A)\mu(B),$$
but this means precisely that the sequence $(\Delta_{T^n})$ converges in the topology of joinings 
to $\mu\otimes\mu$. An elegant application of this characterization was found and used by Ornstein 
in~\cite{ornst1}, where he constructed the first example of a transformation with no square root.

\begin{theo}[Ornstein's criterion for mixing]
The dynamical system $T$ is mixing if and only if $T$ is weakly mixing and there exists some
real number $\theta>0$ such that
\begin{equation}
\label{second-mixing-condition}
\forall A,B\in\A,\ \limsup_{n\to +\infty} \mu(A\cap T^{-n}B)\ \le\ \theta\mu(A)\mu(B).
\end{equation}
\end{theo}
\begin{proof}
The fact that mixing implies weak mixing can be proved by the following argument using joinings: If $T$ is mixing, convergence of $(\Delta_{T^n})$ to $\mu\otimes \mu$ also implies convergence of $(\Delta_{(T\times T)^n})$ (sequence of self-joinings of $T\times T)$ to $\mu\otimes \mu\otimes \mu\otimes \mu$. And this gives that any $T\times T$-invariant set in $X\times X$ has measure zero or one, \textit{i.e.} $T$ is weakly mixing. 
Condition~\eqref{second-mixing-condition} is obviously necessary for $T$ to be mixing. 

Conversely, let us suppose that $T$ is weakly mixing and that~\eqref{second-mixing-condition} holds. Then any cluster point $\lambda$ of the sequence $(\Delta_{T^n})$ in $J(T)$ satisfies $\lambda\le\theta\mu\otimes\mu$, hence is absolutely continuous with respect to $\mu\otimes\mu$. Since $T$ is weakly mixing, $\mu\otimes\mu$ is ergodic, therefore $\lambda=\mu\otimes\mu$. This proves that
$$ \lim_{n\to\infty}\Delta_{T^n}\ =\ \mu\otimes\mu, $$
\textit{i.e.} $T$ is mixing.
\end{proof}


\subsection{2-fold and 3-fold mixing}\label{subsection-2-3-mixing}

We finish this presentation by a long-standing open question
in ergodic theory, and some of its relationships with joining theory.
Recall that the property of (strong) mixing is also called \emph{2-fold mixing},
because it involves only two subsets of the space. A strengthening of this
definition involving 3 (or $k\ge3$) sets gives rise to the property of 
\emph{3-fold mixing} (respectively \emph{$k$-fold mixing}):
\begin{defi}
The dynamical system $T$ is said to be \emph{3-fold mixing} if $\forall A,B,C\in\A$,
$$ \lim_{n,m\to\infty} \mu(A\cap T^{-n}B\cap T^{-(n+m)}C)\ =\ \mu(A)\mu(B)\mu(C). $$
\end{defi}
Again, this definition can easily be translated into the language of joinings: 
$T$ is 3-fold mixing when the sequence $(\Delta_{T^n,T^{n+m}})$ converges in $J^3(T)$
to $\mu\otimes\mu\otimes\mu$ as $n$ and $m$ go to infinity, where $\Delta_{T^n,T^{n+m}}$ 
is the obvious generalization of $\Delta_{T^n}$ to the case of self-joinings of order 3.

Whether 2-fold mixing implies 3-fold mixing was asked long ago by Rohlin \cite{rohli6},
and the question is still open today. However, some important results have been established
around this problem. Some of them could make us think that this implication is
false, by exhibiting counterexamples to the same question but in a generalized context. 
In this category, we must in particular cite Ledrappier's example in \cite{ledra2} of a $\ZZ^2$-action which is   
2-fold but not 3-fold mixing. More recently, it was noticed in \cite{ruedl13} that if
we consider the property of being mixing relatively to a factor $\sigma$-algebra (roughly
speaking, this means that we replace the measure by the conditional expectation with
respect to this $\sigma$-algebra), then we can find examples of 2-fold but not 3-fold mixing.
On the other hand, there are also many results which go in the opposite direction. 
We know that 2-fold mixing implies 3-fold mixing in many cases: For example
when $T$ has singular spectrum (Host, \cite{host1}), when $T$ is rank one (Kalikow, \cite{kalik2})
or even for finite rank systems (Ryzhikov, \cite{ryzhi2}), when $T$ is generated by a Gaussian system
(Leonov, \cite{leono2}; see also Totoki~\cite{totok1}). In fact, in many of these works it is shown that 
a stronger result holds for the respective class of systems under consideration, and
this result concerns joining theory: For finite rank systems, as well as for singular spectrum
systems, it is proved that any ergodic self-joining $\lambda$ of order~3 which has
pairwise independent marginals is necessarily the product measure. As far as Gaussian systems 
are concerned, this property has been established for an important subclass of the zero-entropy Gaussian
processes, called GAG~\cite{leman6}. It is not known whether it can be extended to all zero-entropy Gaussian 
processes. More generally, we may ask the following important question:
\begin{question}
\label{pairwise-independent-question}
Does there exist a zero-entropy dynamical system $T$ and an ergodic joining $\lambda\in J_3(T)$
for which the coordinates are pairwise independent but which is different from $\mu\otimes\mu\otimes\mu$?
\end{question}
Why a negative answer to this question would solve Rohlin's multiple mixing problem 
is explained by the following remarks: 
First, if we had a 2-fold mixing but not 3-fold mixing dynamical system $T$, the sequence 
$(\Delta_{T^n,T^{n+m}})$ would have a cluster point  $\lambda\neq\mu\otimes\mu\otimes\mu$. Since 
$T$ is 2-fold mixing, the coordinates would be pairwise
independent for $\lambda$. And we could find an ergodic joining satisfying those two properties
by taking an ergodic component of $\lambda$. Second, Thouvenot proved that if we could find
a counterexample to Rohlin's question, then we could find one in the zero-entropy category.
This restriction to zero entropy is crucial since we can easily construct examples answering 
Question~\ref{pairwise-independent-question} if we drop it:  
Start from the $(1/2,1/2)$ Bernoulli shift $T$ generated by the i.i.d. process $(\xi_p)_{p\in\ZZ}$,
taking its values in $\{0,1\}$. Take an independent copy $(\xi'_p)_{p\in\ZZ}$ of this process and set
$$ \xi''_p\ \egdef\ \xi_p+\xi'_p\quad \mod 2. $$
Then the three processes $\xi$, $\xi'$ and $\xi''$ have the same distribution, are pairwise independent,
but the self-joining of $T$ of order 3 that we get in this way is not the product measure.

\medskip

Observe also that a negative answer to Question~\ref{pairwise-independent-question} would imply another important result in joining theory: It would show that 
2-fold MSJ implies $k$-fold MSJ for all $k\ge2$. Indeed, any 2-fold MSJ but not 3-fold MSJ would have an ergodic 
joining $\lambda$ of order 3 for which the coordinates would be pairwise independent, but not independent.

\medskip

The purpose of the following is to show that, if ever we could find some system satisfying
the requirements of Question~\ref{pairwise-independent-question}, it must be of a different 
nature than the construction of Ledrappier for an action of $\ZZ^2$, or than
the relative example given in \cite{ruedl13}. Indeed, each of these two constructions leads to
some 
``joining'' $\lambda$, for which we get 3 
pairwise independent coordinates $x_1$, $x_2$ and $x_3$, and for which the global independence 
is denied by the following property: We have some random variable $\xi$ defined on $X$, 
taking its value in $\{0,1\}$ (each value being taken with probability $1/2$), and satisfying  
\begin{equation}
\label{non-independence-equation}
\xi(x_3)\ =\ \xi(x_1) + \xi(x_2)\ \mod 2\quad(\lambda\text{-a.e.}).
\end{equation}

There are two ways of interpreting equation~\eqref{non-independence-equation}:
First, the random variable $\xi$ takes its values in a finite alphabet $A$, and there exists
a function $f:\ A\times A\to A$ such that
$$\xi(x_3)\ =\ f\Bigl(\xi(x_1) , \xi(x_2)\Bigr)\ \quad(\lambda\text{-a.e.}).$$
Second, the random variable $\xi$ takes its values in a compact abelian group $G$, and 
we have
$$ \xi(x_3)\ =\ \xi(x_1) + \xi(x_2)\ \quad(\lambda\text{-a.e.}). $$

We want to show now that there cannot exist an example answering Question~\ref{pairwise-independent-question}
of either of those types, by proving the following two propositions.

\begin{prop}
\label{first-proposition}Suppose that there exist some joining $\lambda\in J_3(T)$, with pairwise independent
coordinates, and some random variable $\xi:\ X\to A$ (finite alphabet) such that 
$$ \xi(x_3)\ =\ f\Bigl(\xi(x_1),\xi(x_2)\Bigr) \quad (\lambda\text{-a.e.}). $$
Then
\begin{itemize}
\item either the factor of $T$ generated by $\xi$ is periodic,
\item or this factor has an entropy at least $\log 2$.
\end{itemize} 
\end{prop}

\begin{prop}
\label{second-proposition}
The dynamical system $T$ is assumed to be ergodic. Suppose that for some joining $\lambda\in J_3(T)$
with pairwise independent marginals, there exists some measurable map $\xi:\ X\to G$ (compact Abelian group) 
which is not $\mu$-a.e. constant, and such that 
$$ \xi(x_3)\ =\ \xi(x_1) + \xi(x_2) \quad (\lambda\mbox{-a.e.}) $$
then $T$ has positive entropy.
\end{prop}

\begin{lemma}
\label{first-lemma}
Let $\xi_1$, $\xi_2$ and $\xi_3$ be 3 random variables taking their values in the 
finite alphabet $A$, identically distributed, pairwise independent, and such that there 
exists a function $f:\ A\times A\to A$ satisfying
$$ \xi_3\ =\ f(\xi_1,\xi_2) \quad (\mbox{a.e.}) $$
Then their common probability law is the uniform law on the set of values which are
taken with positive probability.
\end{lemma}

\begin{proof}
Taking a subset of $A$ instead of $A$ if necessary, we can assume that for all $i\in A$, 
$$p_i:=P(\xi_1=i)>0.$$ 
For all $i$ and $k$ in $A$, there exists a unique $j=j_k(i)$ such that $f(i,j_k(i))=k$.
Moreover, we have
$$ p_k\ =\ \sum_i p_i P(\xi_3=k|\xi_1=i)\ =\ \sum_i p_i p_{j_k(i)}. $$
In other words, $\pi=M\pi$ where $\pi\egdef(p_i)_{i\in A}$ and $M$ is the bistochastic matrix
whose entries are $m_{k,i}:=p_{j_k(i)}$. Taking the limit as $s\to\infty$ in the equality $\pi=M^s\pi$,
we get that $\pi$ has to be the uniform law on $A$. 
\end{proof}

\begin{proof}[Proof of Proposition~\ref{first-proposition}]
From the preceding lemma, we know that $\xi$ is uniformly distributed on the subset of 
values which are taken with positive probability. But for any $m\ge1$, Lemma~\ref{first-lemma} 
also applies to
$$ \xi_0^{m-1}\, :\ x\longmapsto\ \Bigl( \xi(x),\xi(Tx),\ldots,\xi(T^{m-1}x)\Bigr). $$
Let $i_0,\ldots,i_{m-2}$ in $A$ be such that
$$ \mu\Bigl( \xi(x)=i_0, \xi(Tx)=i_1, \ldots,\xi(T^{m-2}x)=i_{m-2} \Bigr)\ >\ 0. $$
Then, the conditional probability law of $\xi(T^{m-1}x)$ knowing the above event is uniform on
the set of values to which it assigns a positive mass. Moreover, we easily check that the cardinal
$a_m$ of the support of this conditional probability law does not depend on the choice of 
$i_0,\ldots,i_{m-2}$, and satisfies $a_m\le a_{m-1}$. Thus,
\begin{itemize}
\item either the sequence $(a_m)$ reaches 1, and in this case the factor generated by $\xi$ is 
periodic;
\item or $a_m$ is always greater than or equal to 2, and then the entropy of the factor generated 
by $\xi$ is at least $\log 2$. 
\end{itemize}
\end{proof}

\begin{lemma}
\label{second-lemma}
Let $\xi_1$, $\xi_2$ and $\xi_3$ be 3 random variables taking their values in the 
compact Abelian group $G$, identically distributed, pairwise independent, and satisfying
$$ \xi_3\ =\ \xi_1+\xi_2 \quad (\mbox{a.e.}) $$
Then their common probability law is the Haar measure on some compact subgroup of $G$.
\end{lemma}

\begin{proof}
Denote by $\nu$ the probability law of  $\xi_1$ on $G$. For $\nu$-almost every $g\in G$, $\nu$ is invariant 
under the addition of $g$, since knowing $\xi_2=g$, both $\xi_1$ and $\xi_3=\xi_1+g$ are distributed according to 
$\nu$. Let us define
$$ H\ \egdef\ \{g\in G:\ \nu\text{ is invariant under the addition of }g\}. $$
This is a closed subgroup of $G$ containing the support of $\nu$. Thus $\nu$ is a probability law
on $H$, invariant under every translation on $H$: $\nu$ is the Haar measure on $H$.
\end{proof}

\begin{proof}[Proof of Proposition~\ref{second-proposition}]
Let us consider $\psi:\ X\to G^{\mathbb{Z}}$ defined by
$$ \psi(x)\ \egdef\ (\xi(T^nx))_{n\in\mathbb{Z}}. $$ 
From Lemma~\ref{second-lemma}, the probability law of $\psi(x)$ is the Haar measure on
some compact subgroup $H$ of $G^\mathbb{Z}$. Moreover, $H$ is invariant by the shift, denoted by
$\sigma$, on $G^\mathbb{Z}$. The restriction of $\sigma$ to $H$ is therefore an ergodic endomorphism of $H$,
which is a factor of $T$ (This is the factor generated by $\xi$). A theorem by Juzvinskii 
\cite{juzvi1} tells us that any such ergodic group endomorphism has to be a $K$-system; hence 
$T$ has positive entropy.
\end{proof}

Observe that the two situations exposed in Propositions~\ref{first-proposition} and~\ref{second-proposition}
can be seen as particular cases of a third one, from which I do not know whether some similar conclusion can 
in general be derived.

\begin{question}
What can be said about a dynamical system $T$ satisfying the following: There exists an ergodic self-joining 
$\lambda\in J_3(T)$, with pairwise independent
marginals, and a map $f:\ X\times X\to X$ such that
$$ x_3\ =\ f(x_1,x_2) \quad(\lambda\text{-a.e.})?$$
\end{question}

\section*{Acknowledgements}

I want to express my special thanks to Jean-Paul Thouvenot and Valery Ryzhikov, from whom I learnt
most of what I know on joinings of dynamical systems; to  Anthony Quas for his help 
concerning Propositions~\ref{first-proposition} and~\ref{second-proposition}; to El Houcein El Abdalaoui, 
\'Elise Janvresse and Yvan Velenik for their support and advice in the preparation of this review.

I am also grateful to the referees for having corrected several mistakes in the first version of this work, and suggested substantial improvements in the text. 

\providecommand{\bysame}{\leavevmode\hbox to3em{\hrulefill}\thinspace}
\providecommand{\MR}{\relax\ifhmode\unskip\space\fi MR }
\providecommand{\MRhref}[2]{%
  \href{http://www.ams.org/mathscinet-getitem?mr=#1}{#2}
}
\providecommand{\href}[2]{#2}

\end{document}